\input amstex
\documentstyle {amsppt}

\topmatter

\title {On abelian subgroups of finite groups} \endtitle

\author {Avinoam Mann} \endauthor

\endtopmatter

\centerline {\bf Abstract}

\medskip

We consider abelain subgroups of small index in finite groups. More generally, we consider subgroups such that the product of their index by the index of their centralizer is small.

\bigskip

\bigskip

In [CD] A.Chermak and A.Delgado introduced a 'measuring argument' for a finite group $G$ acting on a finite group $H$ (all groups in this paper are taken to be finite). A.Goren and M.Herzog generalized it to any permutation action of $G$ [G1, G2, GH]. Both in [CD] and in several later papers (e.g. [Gl]), the case $G = H$, with $G$ acting on itself by conjugation, was singled out, and this case is also the topic of the present paper. We first recall the definitions, adapted for that special case.

For any subgroup $H$ of $G$, we write $m(G,H) := |G:H||G:C_G(H)$, and $\mu(G) := min_{H \le G}m(G,H)$.

{\bf Remark.} Most authors consider the product $|H||C_G(H)|$, but our main interest is in indices, not orders.

A subgroup $H$ such that $m(G,H) = \mu(G)$ is termed a {\it CD-subgroup}.

We state some properties of CD-subgroups.

\proclaim {Theorem 1} Let $H$ and $K$ be CD-subgroups in $G$.

{\bf a.} $HK = KH$, and $C_G(H\cap K) = C_G(H)C_G(K)$.

{\bf b.} $\langle H,K\rangle = HK$, $H\cap K$, and $C_G(H)$ are also CD-subgroups, and $H = C_G(C_G(H))$. Hence the set of CD-subgroups is a lattice, with a unique maximal element $M(G)$ and a unique minimal element $m(G)$. The map $H \rightarrow C_G(H)$ is a 1-1 order reversing map of that lattice onto itself. In particular, $m(G) = C_G(M(G)) = Z(M(G))$ and $M(G) = C_G(m(G))$.

{\bf c.} The CD-subgroups are subnormal in $G$.\endproclaim

Here {\bf a.} and {\bf b.} are proved in [CD], see also [Gl, section 2] for streamlined proofs. {\bf c.} was noted in [Co] and in [BW], and it holds because, by {\bf a.,} a CD-subgroup permutes with its conjugates, and such subgroups are subnormal [Sz].

\proclaim {Corollary 2} $G$ contains a characteristic abelian subgroup of index at most $\mu(G)$.\endproclaim

Indeed, $m(G)$ is such a subgroup.

\proclaim {Corollary 3} If $G/Z(G)$ is a non-abelian simple group, then $G$ and $Z(G)$ are the only CD-subgroups of $G$, and $\mu(G) = |G:Z(G)|$.\endproclaim

This is already stated in [CD] in the case $Z(G) = 1$. To prove it, note that all CD-subgroups are centralizers and contain $Z(G)$, and by assumption no proper subnormal subgroup contains $Z(G)$ properly.

\proclaim {Theorem 4} If $G$ has no non-identity normal abelian subgroups, then its CD-subgroups are exactly its direct factors.\endproclaim

Let $G$ contain an abelian subgroup $A$ of index $n$. Since $A \le C_G(A)$, we have $\mu(G) \le m(G,A) \le n^2$, thus $G$ contains a normal abelian subgroup $N$ of index at most $n^2$. This seems to have been stated first by M.Mozichuk [Mo], and is extended in [PS] to infinite groups. We are going to show that the index can be taken to be strictly less than $n^2$. More generally, we consider nilpotent subgroups.

\proclaim {Theorem 5} If $G$ is not abelian, and $H \neq Z(G)$ is a nilpotent subgroup of $G$, then $G$ contains a normal abelian subgroup of index smaller than $m(G,H)$.\endproclaim

It follows that if $G$ contains an abelian subgroup $A$ of index $n$, it contains a normal abelian subgroup of index at most $\frac{n^2}{p}$, where $p$ is the smallest prime dividing $|G|$. Here it is possible that $G$ has no normal abelian subgroup of smaller index. Our next result determines all $p$-groups with that property. In it we term a group {\it with small abelian subgroups}, if all its abelian subgroups are cyclic ({\it mod} $Z(G)$. These groups are discussed in [He] and [Ma].

\proclaim {Theorem 6} Let $G$ be a non-abelian $p$-groupת in which the minimal index of a normal abelian subgroup is $\frac{\mu(G)}{p}$. Then $G$ satisfies one of the following:

{\bf a.} $G$ contains an abelian maximal subgroup.

{\bf b.} $p$ is odd, and $G/Z(G)$ is the non-abelian group of order $p^3$ and exponent $p$.

{\bf c.} $G$ is a group of class 2 with small abelian subgroups.

\endproclaim

There is a partial converse to this theorem. In the groups satisfying {\bf a.} or {\bf b.}, the minimal index of a normal abelian subgroup is $\frac{\mu(G)}{p}$. Indeed, the non-abelian $p$-groups containing an abelian maximal subgroup are the only $p$-groups in which $\mu(G) = p^2$, and the $p$-groups which do not contain such a maximal subgroup, and in which $|G:Z(G)| = p^3$, are the only $p$-groups in which $\mu(G) = p^3$.
Similarly, it is easily proved that if $\mu(G) = p^4$, then either $|G:Z(G)| = p^4$, or $G$ contains an abelian subgroup of index at most $p^2$. The situation for bigger values of $\mu(G)$ seems to be more complicated.
\bigskip

\centerline {\bf Proofs}

\medskip

{\bf Proof of Theorem 4.} Note that the assumption is equivalent to the assumption that the Fitting subgroup of $G$ is trivial, and so is the soluble radical. Then $G$ also does not contain a non-identity subnormal soluble subgroups. Let $H$ be a CD-subgroup. Then by Theorem 1 $H\cap C_G(H)$ is an abelian subnormal subgroup, hence trivial. Write $D := C_G(H)$. Then $HD = DC_G(D) = C_G(H\cap D) = G$. This shows that $H$ is normal in $G$, and $G = H\times D$.

\medskip

{\bf Proof of Theorem 5.} If $H$ is not a CD-subgroup, the result holds by Corollary 2, so we assume that $H$ is a CD-subgroup. Then it is subnormal, and its normal closure $N$ is also nilpotent. Here

$$\mu(N) \le m(N,H) = |N:H||N:N\cap C_G(H)| = \frac{|G:H|}{|G:N|}\cdot|N:N\cap C_G(H)|$$
$$\le \frac{|G:H|}{|G:N|}\cdot |G:C_G(H)| = \frac{m(G,H)}{|G:N|} = \frac{\mu(G)}{|G:N|}.$$

If $|N:N\cap C_G(H)| < |G:C_G(H)|$, then the second inequality above is strict, and $|G:m(N)| = |N:m(N)||G:N| \le \mu(N)|G:N| < \mu(G)$. Then $m(N)$, which is characteristic in $N$ and normal in $G$, is the normal abelian subgroup that we are looking for.

But if $|N:N\cap C_G(H)| = |G:C_G(H)|$, then $G = NC_G(H)$. Then $N$ is also the normal closure of $H$ in $N$ itself, and in a nilpotent group that means that $N = H$. Thus $G = NC_G(N)$ and $Z(N) = N\cap C_G(N) = Z(G)$ are both CD-subgroups, implying $|G:Z(G)| = \mu(G)$. If $N = G$ then in the nilpotent group $G$ any normal subgroup containing $Z(G)$ with a prime index is normal abelian of index smaller than $|G:Z(G)| = \mu(G)$. Let $N \neq G$. If $N$ is abelian, then $N = Z(N) = Z(G)$, contrary to assumption. Thus the non-abelian nilpotent group $N$ contains a normal abelian subgroup $K$ containing $Z(N)$ properly, so that $|G:K| < |G:Z(N)| = \mu(G)$, and since $G = NC_G(N)$, we have $K \triangleleft G$.

\medskip

{\bf Proof of Theorem 6.} If $|G:M(G)| > p$, then $|G:m(G)| = \frac{\mu(G)}{|G:M(G)|} < \frac{\mu(G)}{p}$. Hence either $M(G) = G$, or $M(G)$ is maximal in $G$. Assume first the latter. If $m(G) = M(G)$, then $M(G)$ is an abelian maximal subgroup. But if $m(G) \neq M(G)$, then there exists a normal subgroup $N \triangleleft G$ such that $m(G) \le N \le M(G)$ and $|N:m(G)| = p$. Since $m(G) = Z(M(G))$, $N$ ia abelian and of index $\frac{|G:m(G)|}{p} = \frac{\mu(G)}{p^2}.$ Thus this case does not occur.
There remains the case $M(G) = G$. Then $m(G) = Z(G)$. Choose a subgroup $N$ as before. Then $N$ is a maximal normal abelian subgroup, hence self-centralizing. Moreover, any maximal normal abelian subgroup has order $p$ ({\it modulo} $Z(G)$. But $N \le Z_2(G)$, and $G'$ centralizes $Z_2(G)$. Thus $G' \le N$ and $cl(G) \le 3$.

If $K$ is another maximal normal abelian subgroup, then $G' \le N\cap K = Z(G)$, and $cl(G) = 2$. Then any maximal abelian subgroup contains $Z(G)$ and is normal, and $G$ has small abelian subgroups.

Finally, if $N$ is the unique maximal normal abelian subgroup, then if $cl(G) = 2$, $G$ still has small abelian subgroups, as above. If $cl(G) = 3$, then $N = Z_2(G)$, because $Z_2(G)$ does not contain elements of order $p^2$ {\it modulo} $Z(G)$. Also $G' \not \le Z(G)$, implying that $G'Z(G)/Z(G) = N/Z(G) = Z(G/Z(G))$, i.e. $G/Z(G)$ is extraspecial. But the only extraspecial groups that are capable of being a central factor group are the dihedral group of order 8, and for an odd prime, the one of order $p^3$ and exponent $p$. In the dihedral case, an element of order 4 in that factor group, together with $Z(G)$, generate an abelian maximal subgroup. This ends the proof.

\bigskip

\centerline {\bf References}

\medskip

BW. B.Brewster-E.Wilcox, Some groups with computable Chermak-Delgado lattices, Bull. Aust. Math. Soc. 86(1) (2012), 29-40.

CD. A.Chermak-A.Delgado, A measuring argument for finite groups, Proc. Amer. Math. Soc. 107 (1989), 907-914.

Co. W.Cooke, Subnormality and the Chermak-Delgado lattice, J. Alg. App. (2020), 7 pages.

G1. A.Goren, A measuring argument for finite permutation groups, Israel J. Math. 145 (2005), 333-339.

G2. A.Goren, Another measuring argument for finite permutation groups, J. Group Th. 10 (2007), 829-840.

GH. A.Goren-M.Herzog, A general measuring argument for finite permutation groups, Proc. Amer. Math. Soc. 137 (2009), 3197-3205.

Gl. G.Glauberman, Centrally large subgroups of finite p-groups, J.Alg. 300 (2006), 480-508.

He. H.Heineken, Gruppen mit kleinen abelschen Untergruppen, Arch. Math. 29 (1977), 20-31.

Ma. A.Mann, Groups with small abelian subgroups, Arch. Math. 50 (1988), 210-213.

Mo. M.Mozichuk, quoted in an early version of [PS].

PS. K.Podosky-B.Szegedy, Bounds in groups with finite abelian coverings or with finite derived groups, J. Gp. Th. 5 (2002), 443-452.

Sz. J.Sz\'ep, Bemerkung zu einem Satz von O.Ore, Publ. Math. Debrecen 3 (1953), 81-82.

\end